\documentclass[12pt]{amsart}

\usepackage{xcolor}	
\usepackage{relsize}
\usepackage{array} 
\usepackage{amsmath}
\usepackage{amssymb}
\usepackage{amsopn}
\usepackage{bm}	
\usepackage{amsthm} 
\usepackage{listings}

\usepackage{graphicx,subcaption}
\usepackage{balance,mathtools}
\usepackage{float} 
\usepackage{tikz-cd} 
\usepackage{hyperref} 
\usepackage{stmaryrd}

\usepackage[backend=biber,hyperref=true,isbn=false,maxnames=7,maxalphanames=4,minalphanames=4,style=alphabetic]{biblatex}
\usepackage{setspace} 

\newcommand{\ve}{\varepsilon}

\newcommand{\mb}{\mathbb}

\newcommand{\RR}{\mb{R}}

\newcommand{\ZZ}{\mb{Z}}

\theoremstyle{definition} 
\newtheorem{defi}{Definition}

\newtheorem{ques}[defi]{Question}
\theoremstyle{plain} 
\newtheorem{thm}[defi]{Theorem}
\newtheorem{cor}[defi]{Corollary}
\newtheorem{lem}[defi]{Lemma}
\newtheorem{conj}[defi]{Conjecture}

\theoremstyle{remark} 
\newtheorem{rem}[defi]{Remark}

\renewbibmacro*{doi+eprint+url}{%
	\printfield{doi}%
	\newunit\newblock%
	\iftoggle{bbx:eprint}{%
		\usebibmacro{eprint}%
	}{}%
	\newunit\newblock%
	\iffieldundef{doi}{%
		\usebibmacro{url+urldate}}%
	{}%
}

\allowdisplaybreaks
\title[Homomesy of tropical $T$-systems for pairs of Dynkin diagrams]{Homomesy of tropical $T$-systems for pairs of simply laced Dynkin diagrams}

\author{Qiyue Chen}

\address{
School of Mathematics \\
	University of Minnesota \\
	Minneapolis\\
	MN 55455\\
	USA
}

\email{chen8448@umn.edu}

\thanks{The author was partially supported by the
Simons Dissertation Fellowship in Mathematics.
}

\date{\today}

\begin{document}

\begin{abstract}
	For a pair of simply laced Dynking diagrams $X$ and $X'$, we prove Chin's homomesy conjecture on the number of  red  and blue mutations in a period of the tropical Zamolodchikov $T$-system associated to $(X, X')$ using Nakanishi's dilogarithm identity.
\end{abstract}
\maketitle
\section{Introduction}

In this paper, we study certain \emph{homomesy}, the phenomenon that some quantity is independent of the orbit under certain automorphism, exhibited in the tropical $T$-systems associated to pairs of simply laced Dynkin diagrams. We fix a pair of simply laced Dynkin diagrams $X_r$ and $X'_{r'}$ with $r$ and $r'$ vertices  and denote $h$ and $h'$ as their Coxeter numbers.

For each $\mathbf{v}=(v,v')\in X_r\times X'_{r'}$, we consider the following \emph{tropical} $T$-system of equations on $\mathbf{t}_{\mathbf{v}}(t)$ with $t\in \mathbb{Z}$:
\begin{equation}\label{eqn:tropical-T}
	\mathbf{t}_{\mathbf{v}}(t)+\mathbf{t}_{\mathbf{v}}(t-2) = \max \left( \sum_{u\sim v} \mathbf{t}_{(u,v')}(t-1), \sum_{v'\sim w'} \mathbf{t}_{(v,w')}(t-1) \right),
\end{equation}
where we denote $u\sim v$ to indicate that $u$ and $v$ are connected by an arrow in $X$, and $v'\sim w'$ to indicate that $v'$ and $w'$ are connected by an arrow in $X'$.

This is the tropical version of the Zamolodchikov $T$-system for pair $(X_r, X'_{r'})$. The famous periodicity theorem of Keller states that:
\begin{thm}[\cite{Keller2013}]
	Any solution to the above system has period dividing $2(h+h')$, i.e.,
	\[
		\mathbf{t}_{\mathbf{v}}(t+2(h+h')) = \mathbf{t}_{\mathbf{v}}(t) \, , \forall \, \mathbf{v}\in X_r\times X'_{r'} , t\in \mathbb{Z}.
	\]
\end{thm}

One can be convinced that if $\mathbf{v}$ and $\mathbf{v'}$ are adjacent to each other in $X_r\times X'_{r'}$, then $\mathbf{t}_{\mathbf{v}}(t_0)$ and $\mathbf{t}_{\mathbf{v'}}(t_0)$ are independent for any $t_0$ in the above system \eqref{eqn:tropical-T}. In fact, the system can be partitioned exactly into two independent subsystems based on parity. Without loss of generality, we color $X_r=I_+\sqcup I_-$ and $X'_{r'}=I_+'\sqcup I_-'$ in a bipartite manner, and introduce the corresponding colorings for the vertices in $X_r\times X'_{r'}$:
\[
	\ve(v,v')=\begin{cases}
		1  & (v,v')\in I_+\times I_+' \sqcup I_-\times I_-', \\
		-1 & (v,v')\in I_+\times I_-' \sqcup I_-\times I_+'.
	\end{cases}
\]
Let $\mathcal{E}_+=\{(v,v',t)\mid \ve(v,v')\cdot (-1)^t = 1\}$. We consider the subsystem defined on $\mathcal{E}_+$, with the following additional initial condition using $\lambda=(\lambda_{\mathbf{v}})\in \RR^{X_r\times X'_{r'}}$:
\begin{equation}
	\begin{split}
		\mathbf{t}_{\mathbf{v}}(0) = \lambda_{\mathbf{v}}  & \impliedby \ve(\mathbf{v})=1,  \\
		\mathbf{t}_{\mathbf{v}}(-1) = \lambda_{\mathbf{v}} & \impliedby \ve(\mathbf{v})=-1.
	\end{split}
\end{equation}

In the classification of Zamolodchikov periodic cluster algebras in \cite{Chin2510.18031}, Chin conjectured the following homomesy in the tropical $T$-systems: If we start with a \emph{generic}
initial labeling $\lambda$ and solve the subsystem on $\mathcal{E}_+$ as above, and look at equations \eqref{eqn:tropical-T} for every
\begin{equation}
	\{(\mathbf{v},t)\in \mathcal{E}_+ \mid 1\le t\le 2(h+h')\}\eqqcolon\mathcal{P}_+,
\end{equation}
we are able to compare whether $\displaystyle\sum_{u\sim v} \mathbf{t}_{(u,v')}(t-1)$ or  $\displaystyle\sum_{v'\sim w'} \mathbf{t}_{(v,w')}(t-1)$ is \emph{strictly} larger. If the former is strictly larger, which means $X$-direction dominates, we call it a \emph{red mutation} at $(\mathbf{v},t)$, and if the latter is strictly larger, which means $X'$-direction dominates, we call it a \emph{blue mutation}.

It is observed and verified for several examples that the total number of red (and blue) mutations over all vertices in one period is independent of the initial labeling (as long as it is generic). More precisely:

\begin{conj}[\cite{Chin2510.18031}, The homomesy conjecture]\label{conj:homomesy}
	Fix a generic labeling. In the $(h+h')rr'$ equations of form \eqref{eqn:tropical-T} for $(\mathbf{v},t)\in \mathcal{P}_+$, there are exactly $hrr'$ red mutations and $h'rr'$ blue mutations for any pair of simply laced Dynkin diagrams $X_r$ and $X'_{r'}$.
\end{conj}
This conjecture has recently been proved in \cite{Chin2610.02548} when $X'_{r'}=A_1$. The main theorem of this paper is:
\begin{thm}\label{thm:main}
	Conjecture \ref{conj:homomesy} is true.
\end{thm}

In \cite{Chin2610.02548}, Chin-Pylyavskyy treated combinatorially the tropical mutations of the system associated to $(X,A_1)$ via laminations on (punctured) polygons, provided an explicit correspondence between the red or blue mutations and lamination curves of certain kind, and counted the number of these curves. They have also constructed a fan structure in the space of initial labelings according to the mutation patterns, and proved that the number of chambers is equal to the number of seeds in the cluster algebra of type $X$.

The proof of Theorem \ref{thm:main} in this paper uses a different approach and relies on a dilogarithm identity established in \cite{Nakanishi2011Jun}. More precisely, we use an idea developed in \cite{Inoue2010Mar,Kuniba2009Dec} which relates $T$-systems to $Y$-systems. Then, using tropicalization technique, we are able to interpret red/blue mutations as whether certain variables in $Y$-system tend to $0/\infty$. And finally, we plug in the $0/\infty$ to the dilogarithm identity of \cite{Nakanishi2011Jun} and deduce the number of red/blue mutations.
This takes us on a detour into the $Y$-systems associated to cluster algebras.

There exist many similarities between the results in \cite{Nakanishi2011Jun} and this paper, which could shed light on an explicit understanding of the homomesy phenomenon similar to the spirit of \cite{Chin2610.02548}. We will discuss the similarities at the end of the paper.

\begin{rem}
	The dilogarithm identity has also been established for non-simply laced Dynkin diagrams, as in \cite{Inoue2013Mar} for type $B_r$ and \cite{Inoue2013MarII} for types $C_r$, $F_4$, and $G_2$.  Similar idea as above applies to Conjecture \ref{conj:homomesy} stated for non-simply laced Dynkin diagrams, as long as we modify the description of (tropical) $T$-systems accordingly. We will leave the precise statement and proof in the future version.
\end{rem}

\begin{rem}
	The full homomesy conjecture in \cite{Chin2510.18031}  considers non-simply laced Dynkin diagrams, and is stated for tropical $T$-systems associated to pairs of Cartan matrices that commute with each other, which aims at covering all the Zamolodchikov periodic cluster algebras. Our Theorem \ref{thm:main} is a special case of this conjecture. We leave the full conjecture for future study.
\end{rem}
The paper is organized as follows: In Section 2, we review the necessary background on cluster algebras, $T$-systems, and $Y$-systems. In Section 3, we explain how to use Nakanishi's dilogarithm identity to prove the above theorem. And in Section 4, we provide further remarks regarding the theorem and its proof.

\subsection*{Acknowledgement}The author would like to thank Ariana Chin for her introduction to this topic at the University's combinatorics seminar and his advisor Gregg Musiker for helpful discussions and useful suggestions.

\section{Background in cluster mutation }

$T$-systems and $Y$-systems are systems of algebraic relations that are widely studied in mathematical physics and integrable systems. In the groundbreaking paper \cite{Fomin2002Apr}, Fomin and Zelevinsky introduced the notion of cluster algebras, under which $T$-systems and $Y$-systems, including the Zamolodchikov system \cite{ZAMOLODCHIKOV1991391}, can be described as exchange relations arising from certain mutation sequences of certain cluster algebras.

\subsection{Cluster Mutations}
A cluster algebra, roughly speaking, is an algebra built from a set of seeds $(B,\mathbf{x},\mathbf{y})$ where $B$ is a skew-symmetric\footnote{We can also study skew-symmetrizable cluster algebras which are beyond the scope of this paper.} matrix associated to a quiver $Q$, $\mathbf{x}$ consists of the cluster variables, and $\mathbf{y}$ consists of the coefficients. In our setting, the $x$-variables corresponds to the variables in the $T$-system, and the $y$-variables correspond to the variables in the $Y$-system. The seeds are related to each other by \emph{mutations} at the vertices of the quiver.
\begin{defi}
	Let $Q$ be a quiver without loops and 2-cycles. Denote by $Q_0$ and $Q_1$ the sets of vertices and arrows of $Q$, respectively. The $B$-matrix of $Q$ is given by $B=(b_{ij})_{i,j\in Q_0}$, where $b_{ij}$ is the number of arrows from $i$ to $j$. If all the arrows are from $j$ to $i$, then $b_{ij}\coloneqq-b_{ji}$.
\end{defi}
To any such quiver, we can associate a cluster algebra. First, we define the \emph{initial cluster} to be $(B,\mathbf{x},\mathbf{y})$ where \begin{itemize}
	\item $B$ is the $B$-matrix associated to the quiver $Q$;
	\item $\mathbf{x}=(x_{i})_{i\in Q_0}$ for some independent indeterminates $x_{i}$ in the field of rational functions.
	\item $\mathbf{y}=(y_{i})_{i\in Q_0}$ for some independent indeterminates $y_{i}$ in the \emph{tropical semifield}.
\end{itemize}

Here, the tropical semifield is the  free abelian group $(\mathbb{P}, \cdot)$ (written multiplicatively) on variables $\{y_i\}_{i\in Q_0}$ and with the \emph{tropical} addition defined by
\begin{equation*}
	\prod_{i\in Q_0} y_i^{\alpha_i}\oplus\prod_{i\in Q_0} y_i^{\beta_i} = \prod_{i\in Q_0} y_i^{\min(\alpha_i, \beta_i)}.
\end{equation*}

\begin{defi}[\cite{Fomin2002Apr,Fomin2007Jan}]
	The seed mutation $\mu_k$ in direction $k$ transforms $(Q,\mathbf{x},\mathbf{y})$ into the seed $\mu_k(Q,\mathbf{x},\mathbf{y}) \coloneqq (Q',\mathbf{x}',\mathbf{y}')$ defined as follows:
	\begin{itemize}

		\item The quiver $Q$ is transformed into a new quiver $Q'$ by the following steps:
		      \begin{enumerate}
			      \item For each subquiver $i\to k\to j$, add a new arrow $i\to j$;
			      \item Reverse all arrows incident to $k$;
			      \item Remove the arrows in a maximal set of pairwise disjoint 2-cycles.
		      \end{enumerate}
		      In terms of the $B$-matrix, the mutation at vertex $k$ transforms $B$ into a new matrix $B'=(b'_{ij})_{i,j\in Q_0}$ by:
		      \[
			      b'_{ij}=\begin{cases}
				      -b_{ij}                                                 & i=k \quad \text{or} \quad j=k          \\
				      b_{ij}+ \operatorname{sign}(b_{ik})\max(b_{ik}b_{kj},0) & i\neq k\quad \text{and} \quad  j\neq k
			      \end{cases}
		      \]
		\item $\mathbf{x}'$ is obtained from $\mathbf{x}$ by replacing one cluster variable by a new one, $\mathbf{x'}=\mathbf{x}\setminus \{x_k\} \cup \{x_k'\}$, and $x_k'$ is defined by the following exchange relation:
		      \[
			      x'_k x_k=
			      \frac{1}{y_k\oplus 1}\left(y_k\prod_{i\to k}x_i+\prod_{k\to i}x_i\right)
		      \]
		\item $\mathbf{y}'=(y'_j)_{j\in Q_0}$ is calculated by:
		      \[
			      y'_j = \begin{cases}
				      \dfrac{1}{y_k}                                                                               & j=k,     \\
				      y_j\displaystyle\prod_{k\to j} y_k (y_k\oplus1)^{-1}\displaystyle\prod_{j\to k}(y_k\oplus 1) & j\neq k.
			      \end{cases}
		      \]

	\end{itemize}
\end{defi}

It can be verified that mutations are involutive. Also, if $b_{jk}=0$, then $\mu_k\mu_j=\mu_j\mu_k$. Moreover, in the study of $T$-systems and $Y$-systems, we either set all the $y$-variables, or $x$-variables, to be equal to 1, and the above mutation rules simplify to:
\begin{itemize}
	\item The exchange relation on $x$-variables:
	      \[x'_k x_k = \prod_{i\to k}x_i + \prod_{k\to i}x_i;\]
	\item The exchange relation on $y$-variables:
	      \[y'_j =\begin{cases}
			      \dfrac{1}{y_k}                                               & j=k,     \\
			      y_j\left(1+y_k^{\operatorname{sign}(b_{jk})}\right)^{b_{jk}} & j\neq k.
		      \end{cases}\]
\end{itemize}

The tropical exchange relation on $x$-variables is obtained by replacing every product with a sum and every sum with a maximum in the ordinary exchange relation:
\[
	\mathbf{t}'_k+\mathbf{t}_k=\max\left(\sum_{i\to k} \mathbf{t}_i, \sum_{k\to i} \mathbf{t}_i\right).
\]

By abuse of notation, the tropical mutations are still denoted by $\mu$ when the context is clear.

\subsection{Recurrent Quivers}
To relate $T$-systems and $Y$-systems to cluster mutations, we start with a bipartite quiver $Q$ with white vertices $Q_+$ and black vertices $Q_-$, such that no arrow connects two vertices of the same color. We define $\mu_+$ and $\mu_-$ as $\prod_{v\in Q_+} \mu_v$ and $\prod_{v\in Q_-} \mu_v$. They are well defined because mutations at non-adjacent vertices commute. We also require $Q$ to be \emph{recurrent}, which means $\mu_+(Q)=\mu_-(Q)=Q^{op}$, where $Q^{op}$ is obtained by reversing all the arrows of $Q$.

For a recurrent quiver $Q$, the associated $T$-system is a set of variables $T_v(t)$ for $v\in Q$ and $t\in \mathbb{Z}$ satisfying the following system of equations:
\[
	T_{v}(t)T_{v}(t-2) = \prod_{u\to v} T_{u}(t-1) + \prod_{v\to w} T_{w}(t-1),
\]
which could also be tropicalized to get the \emph{tropical} $T$-system that we are interested in:
\[
	\mathbf{t}_{v}(t)+\mathbf{t}_{v}(t-2) = \max \left( \sum_{u\to v} \mathbf{t}_{u}(t-1), \sum_{v\to w} \mathbf{t}_{w}(t-1) \right).
\]

The above systems can be partitioned into two independent subsystems, which is similar to the discussion in the introduction. We without loss of generality study the subsystem defined on $Q_+\times 2\ZZ \sqcup (Q_-\times (2\ZZ-1))$.

The cluster mutation provides a solution to the subsystem via the following procedure:
Start by setting
\begin{align*}
	T_v(0)=x_v  & \impliedby v\in Q_+, \\
	T_v(-1)=x_v & \impliedby v\in Q_-,
\end{align*}
then \begin{align*}
	T_v(2t)=(\mu_+ \mu_-)^{t}(\mathbf{x})_v   & \impliedby v\in Q_+, \\
	T_v(2t-1)=(\mu_+ \mu_-)^{t}(\mathbf{x})_v & \impliedby v\in Q_-.
\end{align*}

Following the above procedure of first assigning initial values at $t=0$ or $t=-1$ for vertices, and then applying cluster mutations $\mu_+$ and $\mu_-$, we can also produce the equations of the $Y$-systems. As we will see below in Lemma \ref{lem:square-produce-Za}, $Y$-systems produced in this way generalize the classical Zamolodchikov $Y$-system, which in fact, is the motivation of \cite{Fomin2002Apr} in introducing the concept of general $Y$-system.

\subsection{The Square Product}

In order to study the system \eqref{eqn:tropical-T} for pairs of Dynkin diagrams, \cite{Keller2013} has introduced the \emph{square product quiver} $X_r\square X'_{r'}$. Keep the notations as in the introduction.

For $\mathbf{v}=(v,v')$, write $\mathbf{v}:(++)$ if $v\in I_+$ and $v'\in I'_+$, and so forth for the other combinations of colours. Then,  $X_r\square X'_{r'}$ has vertex set $X_r\times X'_{r'}$, and there is an edge between $(u,v')$ and $(v,v')$ whenever $u$ and $v$ are adjacent in $X$, and between
$(v, v')$ and $(v,w')$ whenever $w$ and $w'$ are adjacent in $
	X'$. Moreover, the edges are oriented in such a way that:
\[
	\begin{tikzcd}[row sep=large, column sep=large]
		++ \arrow[r] & +- \arrow[d] \\
		-+ \arrow[u] & -- \arrow[l]
	\end{tikzcd}
\]
for squares in $X_r\square X'_{r'}$ formed by pairs of adjacent vertices in $X$ and $X'$.

\begin{lem}[\cite{Keller2013}]\label{lem:square-product}
	The quiver $X_r\square X'_{r'}$ defined by the $B$ matrix above is recurrent. Moreover, the associated tropical $T$-system is the system given in \eqref{eqn:tropical-T}.
\end{lem}

In fact, Keller used the square product quiver to establish the periodicity conjecture for the Zamolodchikov $Y$-system \cite{ZAMOLODCHIKOV1991391}.

\begin{lem}[\cite{Keller2013}]\label{lem:square-produce-Za}
	If we start by an initial $y$-seed $(\mathbf{y}_\mathbf{v})$ and assign \begin{equation*}
		\begin{split}
			Y_{\mathbf{v}}(0) = \mathbf{y}_{\mathbf{v}}  & \impliedby \ve(\mathbf{v})=1,  \\
			Y_{\mathbf{v}}(-1) = \mathbf{y}_{\mathbf{v}} & \impliedby \ve(\mathbf{v})=-1,
		\end{split}
	\end{equation*}
	then the variables \emph{defined} by \begin{align*}
		Y_\mathbf{v}(2t)\coloneqq(\mu_+ \mu_-)^{t}(\mathbf{Y})_\mathbf{v}   & \impliedby\ve(\mathbf{v})=1, \\
		Y_\mathbf{v}(2t-1)\coloneqq(\mu_+ \mu_-)^{t}(\mathbf{Y})_\mathbf{v} & \impliedby\ve(\mathbf{v})=-1
	\end{align*}
	satisfy the \emph{Zamolodchikov $Y$-system}
	\begin{equation}\label{eqn:Y-system}
		Y_{v,v'}(t)Y_{v,v'}(t-2)
		=\frac{\displaystyle\prod_{u\sim v} (1+Y_{u,v'}(t-1))}{\displaystyle\prod_{v'\sim w'} (1+Y_{v,w'}(t-1)^{-1})}
	\end{equation}
	for $(\mathbf{v},t)=(v,v',t)\in\mathcal{E}_+$.
\end{lem}

It is based on this cluster algebraic formulation that Keller proves the celebrated periodicity conjecture via a categorification approach through which the periodicity of other cluster algebraic quantities, including the (tropical) $T$-system, can also be established.
\section{Homomesy via Nakanishi's dilogarithm identity}

It is shown in \cite{Nakanishi2011Jun} that the Zamolodchikov $Y$-systems satisfy a remarkable identity involving the Rogers dilogarithm function, and our proof is built on such identity. Recall the Rogers dilogarithm function $L(x)$ defined on $[0,1]$ by
\begin{equation}
	L(x) = -\frac{1}{2} \int_0^x \left(\frac{\log(1-t)}{t}+\frac{\log t}{1-t}\right) \, dt.
\end{equation}

The Rogers dilogarithm function satisfies the following properties:
\begin{itemize}
	\item $L(0) = 0$ and $L(1) = \dfrac{\pi^2}{6}$.
	\item $L$ is continuous on $[0,1]$.
\end{itemize}

\begin{thm}[\cite{Nakanishi2011Jun}, Theorem 2.8]\label{thm:Naka}
	Let $X_r$ and $X'_{r'}$ be simply-laced Dynkin diagrams, with the notations same as above.  For any positive solution of the associated Zamolodchikov $Y$-system \eqref{eqn:Y-system} on $\mathcal{E}_+$, we have
	\begin{equation}\label{eqn:dilog}
		\frac{6}{\pi^2}\sum_{(\mathbf{v},t)\in\mathcal{P}_+}
		L\!\left(\frac{Y_{\mathbf{v}}(t)}{1+Y_{\mathbf{v}}(t)}\right)=hrr'.
	\end{equation}
\end{thm}

\subsection{From $T$-system to $Y$-system}

The ordinary $T$-system for $(\mathbf{v},t)\in \mathcal{E}_+$ assigned to $X_r\square X'_{r'}$ reads
\begin{equation}\label{eqn:T}
	T_{\mathbf{v}}(t)T_{\mathbf{v}}(t-2) = \prod_{u\sim v} T_{(u,v')}(t-1)+ \prod_{v'\sim w'} T_{(v,w')}(t-1) .
\end{equation}
We define two monomials
\begin{align*}
	R_{\mathbf{v}}(t) & =\prod_{u\sim v} T_{(u,v')}(t-1),   \\
	B_{\mathbf{v}}(t) & =\prod_{v'\sim w'} T_{(v,w')}(t-1),
\end{align*}
On $\mathcal{E}_+$, \emph{define}
\begin{equation}\label{eqn:YfromT}
	Y_{\mathbf{v}}(t)=\frac{R_{\mathbf{v}}(t)}{B_{\mathbf{v}}(t)}>0.
\end{equation}

\begin{lem}[\cite{Inoue2010Mar,Kuniba2009Dec}]\label{lem:TtoY}
	For every $(\mathbf{v},t)\in\mathcal{E}_+$, the variables \eqref{eqn:YfromT} satisfy the  Zamolodchikov $Y$-system on $\mathcal{E}_+$:

	\begin{equation}\label{eqn:Y-satisfies-Za}
		Y_{v,v'}(t)Y_{v,v'}(t-2)
		=\frac{\displaystyle\prod_{u\sim v} (1+Y_{u,v'}(t-1))}{\displaystyle\prod_{v'\sim w'} (1+Y_{v,w'}(t-1)^{-1})}.
	\end{equation}

\end{lem}
\begin{proof}
	By \eqref{eqn:YfromT} and the $T$-system applied at the neighbors of $\mathbf{v}$,
	\begin{align*}
		Y_{v,v'}(t)Y_{v,v'}(t-2)
		 & =\frac{\displaystyle\prod_{u\sim v}T_{u,v'}(t-1)T_{u,v'}(t-3)}
		{\displaystyle\prod_{v'\sim w'}T_{v,w'}(t-1)T_{v,w'}(t-3)}                     \\
		 & =\frac{\displaystyle\prod_{u\sim v}\bigl(R_{u,v'}(t-1)+B_{u,v'}(t-1)\bigr)}
		{\displaystyle\prod_{v'\sim w'}\bigl(R_{v,w'}(t-1)+B_{v,w'}(t-1)\bigr)}.
	\end{align*}

	On the other hand, the right-hand side of the equation \eqref{eqn:Y-satisfies-Za} equals the above quotient divided by
	\[
		\frac{\displaystyle\prod_{u\sim v}B_{u,v'}(t-1)}{\displaystyle\prod_{v'\sim w'}R_{v,w'}(t-1)}.
	\]

	This quotient equals $1$, since both its numerator and denominator are equal to
	\[\displaystyle\prod_{v'\sim w'}\displaystyle\prod_{u\sim v}T_{u,w'}(t-2).\qedhere\]
\end{proof}

\subsection{Tropical $T$-system}
As in the introduction, we consider the tropical version of the $T$-system. More precisely, we consider the subsystem defined on $\mathcal{E}_+$, satisfying the equations
\begin{equation*}
	\mathbf{t}_{\mathbf{v}}(t)+\mathbf{t}_{\mathbf{v}}(t-2) = \max \left( \sum_{u\sim v} \mathbf{t}_{(u,v')}(t-1), \sum_{v'\sim w'} \mathbf{t}_{(v,w')}(t-1) \right),
\end{equation*}
with the following  initial condition $\lambda=(\lambda_{\mathbf{v}})\in \RR^{X_r\times X'_{r'}}$:
\begin{equation*}
	\begin{split}
		\mathbf{t}_{\mathbf{v}}(0) = \lambda_{\mathbf{v}}  & \impliedby \ve(\mathbf{v})=1,  \\
		\mathbf{t}_{\mathbf{v}}(-1) = \lambda_{\mathbf{v}} & \impliedby \ve(\mathbf{v})=-1.
	\end{split}
\end{equation*}

Denoting the two entries of the maximum by
\begin{align}\label{eqn:tropical-comparison}
	\mathsf R_{\mathbf{v}}(t) & =\sum_{u\sim v}\mathbf{t}_{u,v'}(t),   \\
	\mathsf B_{\mathbf{v}}(t) & =\sum_{v'\sim w'}\mathbf{t}_{v,w'}(t),
\end{align}
we are interested in the difference
\begin{equation}\label{eqn:D}
	D_{\mathbf{v}}(t)=\mathsf R_{\mathbf{v}}(t)-\mathsf B_{\mathbf{v}}(t).
\end{equation}

\begin{defi}[\cite{Chin2510.18031}]\label{def:colors}
	For the solution of the above tropical system subject to the initial condition, $(\mathbf{v},t)\in\mathcal{P}_+$ is said to give  a \emph{red mutation} if $D_{\mathbf{v}}(t)>0$; a \emph{blue mutation} if $D_{\mathbf{v}}(t)<0$, and a \emph{purple mutation} if $D_{\mathbf{v}}(t)=0$.  A tropical labeling $\lambda$ is \emph{generic} if every $(\mathbf{v},t)\in\mathcal{P}_+$ gives either a red or a blue mutation, i.e., no purple mutation arises in one period.
\end{defi}

For a generic labeling define
$N_{\mathrm{red}}(\lambda)
	=\#\{(\mathbf{v},t)\in\mathcal{P}_H:D_{\mathbf{v}}(t)>0\}$ and
$N_{\mathrm{blue}}(\lambda)
	=\#\{(\mathbf{v},t)\in\mathcal{P}_H:D_{\mathbf{v}}(t)<0\}$.
Then we have $N_{\mathrm{red}}(\lambda)+N_{\mathrm{blue}}(\lambda)=(h+h')rr'$.

Theorem \ref{thm:main} is then the statement that \[
	N_{\mathrm{red}}(\lambda)=hrr',\qquad N_{\mathrm{blue}}(\lambda)=h'rr'.
\]

\subsection{From Ordinary to Tropical}
The following procedure is standard. See for example \cite{Tropicalization}.
For a real parameter $s>0$, specialize the initial labels in the \emph{ordinary} $T$-system by

\begin{equation}\label{eqn:exp-specialization-2}
	\begin{split}
		T^s_{\mathbf{v}}(0) = e^{s\lambda_{\mathbf{v}}}  & \impliedby \ve(\mathbf{v})=1,  \\
		T^s_{\mathbf{v}}(-1) = e^{s\lambda_{\mathbf{v}}} & \impliedby \ve(\mathbf{v})=-1.
	\end{split}
\end{equation}
Let $T^s_{\mathbf{v}}(t)$ and $Y^s_{\mathbf{v}}(t)$ denote the resulting positive solutions of the $T$-system and $Y$-system solved via \eqref{eqn:T} and \eqref{eqn:Y-satisfies-Za}, respectively.

\begin{lem}[The tropical limit]\label{lem:limit}
	For every $(\mathbf{v},t)\in\mathcal{E}_+$, we have:
	\begin{equation}\label{eqn:limits}
		\lim_{s\to+\infty}\frac{\log T^s_{\mathbf{v}}(t)}s=\mathbf{t}_{\mathbf{v}}(t),\qquad
		\lim_{s\to+\infty}\frac{\log Y^s_{\mathbf{v}}(t)}s=D_{\mathbf{v}}(t).
	\end{equation}
\end{lem}
\begin{proof}
	The first identity holds at $t=0$ and $t=-1$ by \eqref{eqn:exp-specialization-2}.  The elementary log-sum-exp limit
	\[
		\lim_{s\to\infty}\frac1s\log\bigl(A_s+B_s\bigr)=\max(\alpha,\beta),
		\quad\text{if}\quad
		\frac{\log A_s}s\to\alpha,\quad \frac{\log B_s}s\to\beta
	\]
	propagates the claim through the recurrence \eqref{eqn:T}, forward and backward, by induction over any finite time interval.  The second identity follows by taking logarithms in \eqref{eqn:YfromT} and using \eqref{eqn:tropical-comparison}--\eqref{eqn:D}.
\end{proof}

\begin{cor}\label{cor:zeroinf}
	For generic $\lambda$ and $(\mathbf{v},t)\in\mathcal{P}_H$,
	\[
		\lim_{s\to\infty}Y_{\mathbf{v}}^s(t)=
		\begin{cases}
			+\infty, & (\mathbf{v},t) \text{ gives a red mutation},  \\
			0,       & (\mathbf{v},t) \text{ gives a blue mutation}.
		\end{cases}
	\]
\end{cor}

\subsection{Proof of Theorem \ref{thm:main}}
Form the one-parameter positive specialization \eqref{eqn:exp-specialization-2}.  Lemma \ref{lem:TtoY} makes the ratios $Y_{\mathbf{v}}^s(t)$ a solution of the  $Y$-system on $\mathcal{E}_+$, so Theorem \ref{thm:Naka} applies for \emph{every} $s>0$:
\begin{equation}\label{eqn:dilog-s}
	\frac{6}{\pi^2}\sum_{(\mathbf{v},t)\in\mathcal{P}_+}
	L\left(\frac{Y_{\mathbf{v}}^s(t)}{1+Y_{\mathbf{v}}^s(t)}\right)=hrr'.
\end{equation}
Because $\lambda$ is generic, Corollary \ref{cor:zeroinf}, together with the fact that $L(0)=0$ and $L(1)=\dfrac{\pi^2}{6}$, show that each summand, after normalization by $6/\pi^2$, tends to $1$ if $(\mathbf{v},t)$ gives a red mutation and $0$ if $(\mathbf{v},t)$ gives a blue mutation.  The sum is finite, so the limit may be passed through the summation sign.  Therefore the left-hand side tends to $N_{\mathrm{red}}(\lambda)$. That $N_{\mathrm{blue}}(\lambda)=h'rr'$ follows from $N_{\mathrm{red}}(\lambda)+N_{\mathrm{blue}}(\lambda)
	=|\mathcal{P}_+|=(h+h')rr'$.\qed

\section{Further remarks}
We remark some similarity in results between \cite{Nakanishi2011Jun} and this paper. In Proposition 3.2, \emph{loc.cit.}, Nakanishi has shown that $h'rr'$ and $hrr'$ are the number of positive and negative tropical $Y$-monomials $[y_{\mathbf{v}}(t)]_\mathrm{T}$ where  $(\mathbf{v},t)\in \mathcal{P}_+$. Moreover, for $0\le t\le h'-1$, $[y_{\mathbf{v}}(t)]_\mathrm{T}$ is a positive monomial when $(\mathbf{v},t)\in \mathcal{E}_+$, and for $-h\le t\le -1$, $[y_{\mathbf{v}}(t)]_\mathrm{T}$ is a negative monomial. The proof of Proposition 3.2 is based on a remarkable \emph{factorization property} of the tropical $Y$-system, which is demonstrated in Example 3.5 \emph{loc.cit.}.

\begin{ques}
	What does the factorization property mean in terms of the tropical $T$-system?
\end{ques}

The proof of Theorem \ref{thm:Naka}, which is Theorem 2.8 in \cite{Nakanishi2011Jun}, is based on this Proposition 3.2, as well as an application of the $0/\infty$-technique in its Section 4.2 which is similar to our proof strategy as in Corollary \ref{cor:zeroinf} and below. Given these similarities, we would like to ask:

\begin{ques}
	What is the relationship between those positive or negative tropical $Y$-monomials $[y_{\mathbf{v}}(t)]_\mathrm{T}$ and the set of blue or red mutations?
\end{ques}

On the other hand, as is mentioned in the introduction, \cite{Chin2610.02548} provided an explicit combinatorial description of the red and blue mutations for systems associated to $(X,A_1)$ via lamination curves on (punctured) polygons, and proved the conjecture via counting these curves. For questions related to finite type Dynkin diagrams, the polygonal understanding continuously proves to be useful. On the other hand, Keller's proof of periodicity for systems associated to general $(X,X')$ is heavily based on representation theory and categorification, which motivates us to ask:

\begin{ques}
	How to understand the red and blue mutations for pairs of Dynkin diagrams in a more explicit combinatorial, representation-theoretical or categorical way?
\end{ques}

\printbibliography
\end{document}